\documentclass[11pt, a4paper]{article}
\usepackage[cp1251]{inputenc}
\usepackage[russian]{babel}
\usepackage{amsmath} \usepackage{euscript}
\usepackage{mymatrix2}%\let\mymatrixII\mymatrix
\usepackage{longtable}
\usepackage{mathrsfs}
\usepackage{amssymb}
\begin{document}

\newcounter{bnomer} \newcounter{snomer}
\newcounter{bsnomer}
\setcounter{bnomer}{0}
\renewcommand{\thesnomer}{\thebnomer.\arabic{snomer}}
\renewcommand{\thebsnomer}{\thebnomer.\arabic{bsnomer}}
\renewcommand{\refname}{\begin{center}\large{\textbf{References}}\end{center}}

\setcounter{MaxMatrixCols}{14}

\newcommand{\sect}[1]{%
\setcounter{snomer}{0}\setcounter{bsnomer}{0}
\refstepcounter{bnomer}
\par\bigskip\begin{center}\large{\textbf{\arabic{bnomer}. {#1}}}\end{center}}
\newcommand{\sst}{%
\refstepcounter{bsnomer}
\par\bigskip\textbf{\arabic{bnomer}.\arabic{bsnomer}. }}
\newcommand{\defi}[1]{%
\refstepcounter{snomer}
\par\medskip\textbf{Definition \arabic{bnomer}.\arabic{snomer}. }{#1}\par\medskip}
\newcommand{\theo}[2]{%
\refstepcounter{snomer}
\par\textbf{Теорема \arabic{bnomer}.\arabic{snomer}. }{#2} {\emph{#1}}\hspace{\fill}$\square$\par}
\newcommand{\mtheop}[2]{%
\refstepcounter{snomer}
\par\textbf{Theorem \arabic{bnomer}.\arabic{snomer}. }{\emph{#1}}
\par\textsc{Proof}. {#2}\hspace{\fill}$\square$\par}
\newcommand{\mcorop}[2]{%
\refstepcounter{snomer}
\par\textbf{Corollary \arabic{bnomer}.\arabic{snomer}. }{\emph{#1}}
\par\textsc{Proof}. {#2}\hspace{\fill}$\square$\par}
\newcommand{\mtheo}[1]{%
\refstepcounter{snomer}
\par\medskip\textbf{Theorem \arabic{bnomer}.\arabic{snomer}. }{\emph{#1}}\par\medskip}
\newcommand{\mlemm}[1]{%
\refstepcounter{snomer}
\par\medskip\textbf{Lemma \arabic{bnomer}.\arabic{snomer}. }{\emph{#1}}\par\medskip}
\newcommand{\mprop}[1]{%
\refstepcounter{snomer}
\par\medskip\textbf{Proposition \arabic{bnomer}.\arabic{snomer}. }{\emph{#1}}\par\medskip}
\newcommand{\theobp}[2]{%
\refstepcounter{snomer}
\par\textbf{Теорема \arabic{bnomer}.\arabic{snomer}. }{#2} {\emph{#1}}\par}
\newcommand{\theop}[2]{%
\refstepcounter{snomer}
\par\textbf{Theorem \arabic{bnomer}.\arabic{snomer}. }{\emph{#1}}
\par\textsc{Proof}. {#2}\hspace{\fill}$\square$\par}
\newcommand{\theosp}[2]{%
\refstepcounter{snomer}
\par\textbf{Теорема \arabic{bnomer}.\arabic{snomer}. }{\emph{#1}}
\par\textbf{Схема доказательства}. {#2}\hspace{\fill}$\square$\par}
\newcommand{\exam}[1]{%
\refstepcounter{snomer}
\par\medskip\textbf{Example \arabic{bnomer}.\arabic{snomer}. }{#1}\par\medskip}
\newcommand{\deno}[1]{%
\refstepcounter{snomer}
\par\textbf{Definition \arabic{bnomer}.\arabic{snomer}. }{#1}\par}
\newcommand{\post}[1]{%
\refstepcounter{snomer}
\par\textbf{Предложение \arabic{bnomer}.\arabic{snomer}. }{\emph{#1}}\hspace{\fill}$\square$\par}
\newcommand{\postp}[2]{%
\refstepcounter{snomer}
\par\medskip\textbf{Proposition \arabic{bnomer}.\arabic{snomer}. }{\emph{#1}}%
\ifhmode\par\fi\textsc{Proof}. {#2}\hspace{\fill}$\square$\par\medskip}
\newcommand{\lemm}[1]{%
\refstepcounter{snomer}
\par\textbf{Lemma \arabic{bnomer}.\arabic{snomer}. }{\emph{#1}}\hspace{\fill}$\square$\par}
\newcommand{\lemmp}[2]{%
\refstepcounter{snomer}
\par\medskip\textbf{Lemma \arabic{bnomer}.\arabic{snomer}. }{\emph{#1}}
\par\textsc{Proof}. {#2}\hspace{\fill}$\square$\par\medskip}
\newcommand{\coro}[1]{%
\refstepcounter{snomer}
\par\textbf{Следствие \arabic{bnomer}.\arabic{snomer}. }{\emph{#1}}\hspace{\fill}$\square$\par}
\newcommand{\mcoro}[1]{%
\refstepcounter{snomer}
\par\textbf{Corollary \arabic{bnomer}.\arabic{snomer}. }{\emph{#1}}\par\medskip}
\newcommand{\corop}[2]{%
\refstepcounter{snomer}
\par\textbf{Следствие \arabic{bnomer}.\arabic{snomer}. }{\emph{#1}}
\par\textsc{Proof}. {#2}\hspace{\fill}$\square$\par}
\newcommand{\nota}[1]{%
\refstepcounter{snomer}
\par\medskip\textbf{Remark \arabic{bnomer}.\arabic{snomer}. }{#1}\par\medskip}
\newcommand{\propp}[2]{%
\refstepcounter{snomer}
\par\medskip\textbf{Proposition \arabic{bnomer}.\arabic{snomer}. }{\emph{#1}}
\par\textsc{Proof}. {#2}\hspace{\fill}$\square$\par\medskip}
\newcommand{\hypo}[1]{%
\refstepcounter{snomer}
\par\medskip\textbf{Conjecture \arabic{bnomer}.\arabic{snomer}. }{\emph{#1}}\par\medskip}
\newcommand{\prop}[1]{%
\refstepcounter{snomer}
\par\textbf{Proposition \arabic{bnomer}.\arabic{snomer}. }{\emph{#1}}\hspace{\fill}$\square$\par}

\newcommand{\Ind}[3]{%
\mathrm{Ind}_{#1}^{#2}{#3}}
\newcommand{\Res}[3]{%
\mathrm{Res}_{#1}^{#2}{#3}}
\newcommand{\epsi}{\epsilon}
\newcommand{\tri}{\triangleleft}
\newcommand{\Supp}[1]{%
\mathrm{Supp}(#1)}

\newcommand{\reg}{\mathrm{reg}}
\newcommand{\empr}[2]{[-{#1},{#1}]\times[-{#2},{#2}]}
\newcommand{\sreg}{\mathrm{sreg}}
\newcommand{\codim}{\mathrm{codim}\,}
\newcommand{\chara}{\mathrm{char}\,}
\newcommand{\rk}{\mathrm{rk}\,}
\newcommand{\chr}{\mathrm{ch}\,}
\newcommand{\id}{\mathrm{id}}
\newcommand{\Ad}{\mathrm{Ad}}
\newcommand{\col}{\mathrm{col}}
\newcommand{\row}{\mathrm{row}}
\newcommand{\low}{\mathrm{low}}
\newcommand{\pho}{\hphantom{\quad}\vphantom{\mid}}
\newcommand{\fho}[1]{\vphantom{\mid}\setbox0\hbox{00}\hbox to \wd0{\hss\ensuremath{#1}\hss}}
\newcommand{\wt}{\widetilde}
\newcommand{\wh}{\widehat}
\newcommand{\ad}[1]{\mathrm{ad}_{#1}}
\newcommand{\tr}{\mathrm{tr}\,}
\newcommand{\GL}{\mathrm{GL}}
\newcommand{\SL}{\mathrm{SL}}
\newcommand{\SO}{\mathrm{SO}}
\newcommand{\Sp}{\mathrm{Sp}}
\newcommand{\Mat}{\mathrm{Mat}}
\newcommand{\Stab}{\mathrm{Stab}}

\newcommand{\vfi}{\varphi}
\newcommand{\teta}{\vartheta}
\newcommand{\Bfi}{\Phi}
\newcommand{\Fp}{\mathbb{F}}
\newcommand{\Rp}{\mathbb{R}}
\newcommand{\Zp}{\mathbb{Z}}
\newcommand{\Cp}{\mathbb{C}}
\newcommand{\ut}{\mathfrak{u}}
\newcommand{\at}{\mathfrak{a}}
\newcommand{\nt}{\mathfrak{n}}
\newcommand{\mt}{\mathfrak{m}}
\newcommand{\htt}{\mathfrak{h}}
\newcommand{\spt}{\mathfrak{sp}}
\newcommand{\rt}{\mathfrak{r}}
\newcommand{\rad}{\mathfrak{rad}}
\newcommand{\bt}{\mathfrak{b}}
\newcommand{\gt}{\mathfrak{g}}
\newcommand{\vt}{\mathfrak{v}}
\newcommand{\pt}{\mathfrak{p}}
\newcommand{\Xt}{\mathfrak{X}}
\newcommand{\Po}{\mathcal{P}}
\newcommand{\Uo}{\EuScript{U}}
\newcommand{\Fo}{\EuScript{F}}
\newcommand{\Do}{\EuScript{D}}
\newcommand{\Eo}{\EuScript{E}}
\newcommand{\Iu}{\mathcal{I}}
\newcommand{\Mo}{\mathcal{M}}
\newcommand{\Nu}{\mathcal{N}}
\newcommand{\Ro}{\mathcal{R}}
\newcommand{\Co}{\mathcal{C}}
\newcommand{\Lo}{\mathcal{L}}
\newcommand{\Ou}{\mathcal{O}}
\newcommand{\Uu}{\mathcal{U}}
\newcommand{\Au}{\mathcal{A}}
\newcommand{\Vu}{\mathcal{V}}
\newcommand{\Bu}{\mathcal{B}}
\newcommand{\Sy}{\mathcal{Z}}
\newcommand{\Sb}{\mathcal{F}}
\newcommand{\Gr}{\mathcal{G}}
\newcommand{\rtc}[1]{C_{#1}^{\mathrm{red}}}

\author{Mikhail V. Ignatyev\and Aleksandr A. Shevchenko}

\date{Samara State University, Chair of algebra and geometry\\\texttt{mihail.ignatev@gmail.com}\\\texttt{shevchenko.alexander.1618@gmail.com}}
\title{\Large{On tangent cones to Schubert varieties in type $D_n$}\mbox{$\vphantom{1}$}\footnotetext{The authors were partially supported by RFBR grants no. 14--01--31052 and 14--01--97017. The first author was partially supported by the Dynasty Foundation, by Max Planck Institute for Mathematics and by the Ministry of Science and Education of the Russian
Federation.}} \maketitle

\sect{Introduction and the main results}

\sst Let $G$ be a complex reductive algebraic group, $T$ a maximal torus in~$G$, $B$ a Borel subgroup in~$G$ containing $T$, and $U$ the unipotent radical of $B$. Let $\Phi$ be the root system of $G$ with respect to~$T$, $\Phi^+$ the set of positive roots with respect to $B$, $\Delta$ the set of simple roots, and $W$ the Weyl group of $\Phi$ (see \cite{Bourbaki}, \cite{Humphreys} and \cite{Humpreys2} for basic facts about algebraic groups and root systems).

Denote by $\Fo=G/B$ the flag variety and by $X_w\subseteq\Fo$ the Schubert subvariety corresponding to an element $w$ of the Weyl group $W$. Denote by $\Ou=\Ou_{p,X_w}$ the local ring at the point $p=eB\in X_w$. Let $\mt$ be the maximal ideal of~$\Ou$. The sequence of ideals $$\Ou\supseteq\mt\supseteq\mt^2\supseteq\ldots$$ is a filtration on $\Ou$. We define $R$ to be the graded algebra $$R=\mathrm{gr}\,\Ou=\bigoplus_{i\geq0}\mt^i/\mt^{i+1}.$$ By definition, the \emph{tangent cone} $C_w$ to the Schubert variety $X_w$ at the point $p$ is the spectrum of~$R$: $C_w=\mathrm{Spec}\,R$. Obviously, $C_w$ is a subscheme of the tangent space $T_pX_w\subseteq T_p\Fo$. A hard problem in studying geometry of $X_w$ is to describe $C_w$ \cite[Chapter 7]{BilleyLakshmibai}.

In 2011, D.Yu. Eliseev and A.N. Panov computed tangent cones $C_w$ for all $w\in W$ in the case $G=\mathrm{SL}_n(\mathbb{C})$, $n\leq5$ \cite{EliseevPanov}. Using their computations, A.N. Panov formulated the following Conjecture.

\hypo{\textup{(A.N. Panov, 2011)} Let $w_1$\textup{,} $w_2$ be \label{mconj}involutions\textup{,} i.e.\textup{,} $w_1^2=w_2^2=\id$. If $w_1\neq w_2$\textup{,} then $C_{w_1}\neq C_{w_2}$ as subschemes of $T_p\Fo$.}

One can easily check that it is enough to prove the Conjecture for irreducible root systems (see Remark~\ref{nota:irred} below). In 2013, D.Yu. Eliseev and the first author proved this Conjecture in types $A_n$, $F_4$ and $G_2$ \cite{EliseevIgnatyev}. In \cite{BochkarevIgnatyevShevchenko}, M.A. Bochkarev and the authors proved the Conjecture in types $B_n$ and $C_n$. In this paper, we prove that the Conjecture is true if $\Phi$ is of type $D_n$ and $w_1$, $w_2$ are basic involutions (see Definition~\ref{defi:basic_involution}). Precisely, our first main result is as follows.

\mtheo{Assume that every irreducible component of $\Phi$ is of type $D_n$ \textup, $n\geq4$. Let $w_1$\textup,~$w_2$ be basic involutions in the Weyl group of $\Phi$ and $w_1\neq w_2$. Then the tangent cones $C_{w_1}$ and~$C_{w_2}$ do~not coincide as subschemes of $T_p\Fo$.\label{mtheo:non_red}}

Note that the similar question for other involutions in $D_n$ and for the root systems $E_6$, $E_7$, $E_8$ remains open.

Now, let $\Au$ be the symmetric algebra of the vector space $\mt/\mt^2$, or, equivalently, the algebra of regular functions on the tangent space $T_pX_w$. Since $R$ is generated as $\Cp$-algebra by $\mt/\mt^2$, it is a~quotient ring $R=\Au/I$. By definition, the \emph{reduced tangent cone} $C_w^{\mathrm{red}}$ to $X_w$ at the point $p$ is the common zero locus in $T_pX_w$ of the polynomials $f\in I\subseteq\Au$. Clearly, if $\rtc{w_1}\neq\rtc{w_2}$, then $C_{w_1}\neq C_{w_2}$. Our second main result is as follows.

\mtheo{Assume that every irreducible component of $\Phi$ is of type $D_n$\textup, $n\geq4$. Let $w_1$\textup, $w_2$ be basic involutions in the Weyl group of $\Phi$ and $w_1\neq w_2$. Then the reduced tangent cones $\rtc{w_1}$ and $\rtc{w_2}$ do not coincide as subvarieties of $T_p\Fo$.\label{mtheo:red}}

In \cite{BochkarevIgnatyevShevchenko}, the similar result was obtained by M.A. Bochkarev for root systems of types $A_n$ and $C_n$. Our proof for $D_n$ is based on the similar idea. Note that the similar question for other involutions in $D_n$ and for other root systems remains open.

The paper is organized as follows. In the next Subsection, we introduce the main technical tool used in the proof of Theorem~\ref{mtheo:non_red}. Namely, to each element $w\in W$ one can assign a polynomial $d_w$ in the algebra of regular functions on the Lie algebra of the maximal torus $T$. These polynomials are called Kostant--Kumar polynomials \cite{KostantKumar1}, \cite{KostantKumar2}, \cite{Kumar}, \cite{Billey}. In \cite{Kumar} S. Kumar showed that if $w_1$ and~$w_2$ are arbitrary elements of~$W$ and $d_{w_1}\neq d_{w_2}$, then $C_{w_1}\neq C_{w_2}$. We give three equivalent definitions of Kostant--Kumar polynomials and formulate their properties needed for the sequel.

In Section~\ref{sect:non_red} we prove that if all irreducible components of $\Phi$ are of type $D_n$ and $w_1$, $w_2$ are distinct basic involutions in $W$, then $d_{w_1}\neq d_{w_2}$, see Proposition~\ref{prop:non_red_eq_C_1}. This implies that $C_{w_1}\neq C_{w_2}$ and proves Theorem~\ref{mtheo:non_red}. The proofs of Conjecture~\ref{mconj} for $A_n$, $F_4$, $G_2$, $B_n$ and $C_n$ presented in \cite{EliseevIgnatyev} and \cite{BochkarevIgnatyevShevchenko} are based on the similar argument.

Section~\ref{sect:red} contains the proof of Theorem~\ref{mtheo:red}. Namely, in Subsection~\ref{sst:red_coadjoint} we describe connections of the geometry of tangent cones with the geometry of coadjoint $B$-orbits. Using these connections, in Subsection~\ref{sst:red_D_n} we proof the result.

Of course, Theorem~\ref{mtheo:non_red} is a corollary of Theorem~\ref{mtheo:red}. Nevertheless, we give in Section~\ref{sect:non_red} an independent proof of the first Theorem based on computation of Kostant--Kumar polynomials. The reason is that we hope to prove Theorem~\ref{mtheo:non_red} for all involutions in $D_n$ using the same technique. At the contrary, there is no chance to prove Theorem~\ref{mtheo:red} for non-basic involutions in $D_n$ using arguments similar to presented in Section~\ref{sect:red}, see Remark~\ref{nota:why_basic_red} (ii) for the details.

\medskip\textsc{Acknowledgements}. The authors were supported by RFBR grants no. 14--01--31052 and\break 14--01--97017; RFBR is gratefully acknowledged. Mikhail Ignatyev aknowledges support from the Dy\-nas\-ty Foundation and from the Ministry of Science and Education of the Russian Federation. The work was completed during the stay of Mikhail Ignatyev at Max Planck Institute for Mathematics in September--October 2014. Mikhail Ignatyev thanks MPIM for the hospitality and for the financial support.

\sst Let $w$ be an element of the Weyl group $W$. Here we give precise definition of the Kostant--Kumar polynomial $d_w$, explain how to compute it in combinatorial terms, and show that it depends only on the scheme structure of $C_w$.

The torus $T$ acts on the Schubert variety $X_w$ by left multiplications (or, equivalently, by con\-ju\-ga\-tions). The point $p$ is invariant under this action, hence there is the structure of a $T$-module on the local ring $\Ou$. The action of $T$ on $\Ou$ preserves the filtration by powers of the ideal $\mt$, so we obtain the structure of a $T$-module on the algebra $R=\mathrm{gr}\,\Ou$. By~\cite[Theorem~2.2]{Kumar}, $R$ can be decomposed into a direct sum of its finite-dimensional weight subspaces: $$R=\bigoplus_{\lambda\in\Xt(T)}R_{\lambda}.$$
Here $\htt$ is the Lie algebra of the torus $T$, $\Xt(T)\subseteq\htt^*$ is the character lattice of $T$ and $R_{\lambda}=\{f\in R\mid t.f=\lambda(t)f\}$ is the weight subspace of weight $\lambda$. Let $\Lambda$ be the $\Zp$-module consisting of all (possibly infinite) $\Zp$-linear combinations of linearly independent elements $e^{\lambda}$, $\lambda\in\Xt(T)$. The \emph{formal character} of~$R$ is an element of $\Lambda$ of the form $$\chr R=\sum_{\lambda\in\Xt(T)}m_{\lambda}e^{\lambda},$$ where $m_{\lambda}=\dim R_{\lambda}$.

Now, pick an element $a=\sum_{\lambda\in\Xt(T)}n_{\lambda}e^{\lambda}\in\Lambda$. Assume that there are finitely many $\lambda\in\Xt(T)$ such that $n_{\lambda}\neq0$. Given $k\geq0$, one can define the polynomial $$[a]_k=\sum_{\lambda\in\Xt(T)}n_{\lambda}\cdot\dfrac{\lambda^k}{k!}\in S=\Cp[\htt].$$ Denote $[a]=[a]_{k_0}$, where $k_0$ is minimal among all non-negative numbers $k$ such that $[a]_k\neq0$. For instance, if $a=1-e^{\lambda}$, then $[a]_0=0$ and $[a]=[a]_1=-\lambda$ (here we denote $1=e^0$).

Let $A$ be the submodule of $\Lambda$ consisting of all finite linear combinations. It is a commutative ring with respect to the multiplication $e^{\lambda}\cdot e^{\mu}=e^{\lambda+\mu}$. In fact, it is just the group ring of $\Xt(T)$. Denote the field of fractions of the ring $A$ by $Q\subseteq\Lambda$. To each element of $Q$ of the form $q=a/b$, $a$,~$b\in A$, one can assign the element $$[q]=\dfrac{[a]}{[b]}\in\Cp(\htt)$$ of the field of rational functions on $\htt$. Note that this element is well-defined~\cite{Kumar}.

There exists an involution $q\mapsto q^*$ on $Q$ defined by $$e^{\lambda}\mapsto (e^{\lambda})^*=e^{-\lambda}.$$ It turns out \cite[Theorem 2.2]{Kumar} that the character $\chr R$ belongs to $Q$, hence $(\chr R)^*\in Q$, too. Finally, we put $$c_w=[(\chr R)^*],\ d_w=(-1)^{l(w)}\cdot c_w\cdot\prod_{\alpha\in\Phi^+}\alpha.$$ Here $l(w)$ is the length of $w$ in the Weyl group $W$ with respect to the set of simple roots $\Delta$. Evidently, $c_w$ and $d_w$ belong to $\Cp(\htt)$; in fact, $d_w$ is a polynomial, i.e., it belongs to the algebra $S=\Cp[\htt]$ of regular functions on $\htt$, see \cite{KostantKumar2} and \cite[Theorem 7.2.6]{BilleyLakshmibai}. \defi{Let $w$ be an element of the Weyl group $W$. The polynomial $d_w\in S$ is called the \emph{Kostant--Kumar polynomial} associated with $w$.}

It follows from the definition that $c_w$ and $d_w$ depend only on the canonical structure of a $T$-module on the algebra $R$ of regular functions on the tangent cone $C_w$. Thus, to prove that the tangent cones corresponding to elements $w_1$, $w_2$ of the Weyl group are distinct, it is enough to check that ${c_{w_1}\neq c_{w_2}}$, or, equivalently, $d_{w_1}\neq d_{w_2}$.

On the other hand, there is a purely combinatorial description of Kostant--Kumar polynomials. To give this description, we need some more notation. Let $w$, $v$ be elements of~$W$. Fix a reduced decomposition of the element $w=s_{i_1}\ldots s_{i_l}$. (Here $\alpha_1,\ldots,\alpha_n\in\Delta$ are simple roots and $s_i$ is the simple reflection corresponding to $\alpha_i$.) Put
\begin{equation*}
c_{w,v}=(-1)^{l(w)}\cdot\sum\dfrac{1}{s_{i_1}^{\epsi_1}\alpha_{i_1}}\cdot\dfrac{1}{s_{i_1}^{\epsi_1}
s_{i_2}^{\epsi_2}\alpha_{i_2}}\cdot\ldots\cdot\dfrac{1}{s_{i_1}^{\epsi_1}\ldots s_{i_l}^{\epsi_l}\alpha_{i_l}},
\end{equation*}
where the sum is taken over all sequences $(\epsi_1,\ldots,\epsi_l)$ of zeroes and units such that ${s_{i_1}^{\epsi_1}\ldots s_{i_l}^{\epsi_l}=v}$. Actually, the element $c_{w,v}\in\Cp(\htt)$ depends only on $w$ and $v$, not on the choice of a reduced de\-com\-po\-si\-tion of $w$ \cite[Section 3]{Kumar}.

\exam{Let $\Phi=A_n$. Put $w=s_1s_2s_1$. To compute $c_{w,\id}$, we should take the sum over two sequences, $(0,0,0)$ and $(1,0,1)$. Hence
\begin{equation*}
c_{w,\id}=(-1)^3\cdot\left(\dfrac{1}{\alpha_1\alpha_2\alpha_1}+\dfrac{1}{-\alpha_1(\alpha_1+\alpha_2)\alpha_1}\right)
=\dfrac{1}{\alpha_1\alpha_2(\alpha_1+\alpha_2)}.
\end{equation*}}

A remarkable fact is that $c_{w,\id}=c_w$, hence to prove that the tangent cones to Schubert varieties do not coincide as subschemes, we need only combinatorics of the Weyl group. Note also that for classical Weyl groups, elements $c_{w,v}$ are closely related to Schubert polynomials \cite{Billey}.

Finally, we will present an original definition of elements $c_{w,v}$ using so-called nil-Hecke ring (see \cite{Kumar} and~\cite[Section 7.1]{BilleyLakshmibai}). The group $W$ naturally acts on $\Cp(\htt)$ by automorphisms. Denote by $Q_W$ the vector space over $\Cp(\htt)$ with basis $\{\delta_w,\ w\in W\}$. It is a ring with respect to the multiplication $$f\delta_v\cdot g\delta_w=fv(g)\delta_{vw}.$$
This ring is called the \emph{nil-Hecke ring}. To each $i$ from 1 to $n$ put $$x_i=\alpha_i^{-1}(\delta_{s_i}-\delta_{\id}).$$
Let $w\in W$ and $w=s_{i_1}\ldots s_{i_l}$ be a reduced decomposition of $w$. Then the element $$x_w=x_{i_1}\ldots x_{i_l}$$ does not depend on the choice of a reduced decomposition of $w$ \cite[Proposition 2.1]{KostantKumar1}.

Moreover, it turns out that $\{x_w,\ w\in W\}$ is a $\Cp(\htt)$-basis of~$Q_W$ \cite[Proposition 2.2]{KostantKumar1}, and
\begin{equation*}
x_w=\sum\nolimits_{v\in W}c_{w,v}\delta_v.\\
\end{equation*}
Actually, if $w,v\in W$, then
\begin{equation}
\begin{split}
&\text{a) }x_v\cdot x_w=\begin{cases}x_{vw},&\text{ if }l(vw)=l(v)+l(w),\\
0,&\text{ otherwise},
\end{cases}\\
&\text{b) }c_{w,v}=-v(\alpha_i)^{-1}(c_{ws_i,v}+c_{ws_i,vs_i}),\text{ if }l(ws_i)=l(w)-1,\\
&\text{c) }c_{w,v}=\alpha_i^{-1}(s_i(c_{s_iw,s_iv})-c_{s_iw,v}),\text{ if }l(s_iw)=l(w)-1.\\
\end{split}\label{formula_x_v_x_w}
\end{equation}
The first property is proved in \cite[Proposition 2.2]{KostantKumar1}. The second and the third properties follow immediately from the first one and the definitions (see also the proof of \cite[Corollary 3.2]{Kumar}).

\nota{i) Suppose~$\Phi$ is a union of its\label{nota:irred} subsystems $\Phi_1$ and $\Phi_2$ contained in mutually orthogonal subspaces. Let $W_1$, $W_2$ be the Weyl groups of $\Phi_1$, $\Phi_2$ respectively, so $W=W_1\times W_2$. Denote $\Delta_1=\Delta\cap\Phi_1=\{\alpha_1,\ldots,\alpha_r\}$ and $\Delta_2=\Delta\cap\Phi_2=\{\beta_1,\ldots,\beta_s\}$, then $$\Cp[\htt]\cong\Cp[\alpha_1,\ldots,\alpha_r,\beta_1,\ldots,\beta_s].$$
Given $v\in W_1$, denote by $d_v^1$ its Kostant--Kumar polynomial. We can consider $d_v^1$ as an element of~$\Cp(\htt)$ depending only on $\alpha_1,\ldots,\alpha_r$. We define $c_v^1\in\Cp(\htt)$ by the similar way. Given $v\in W_2$, we define $d_v^2\in\Cp[\htt]$ and $c_v^2\in\Cp(\htt)$; they depend only on $\beta_1,\ldots,\beta_s$. Let $w\in W$, $w_1\in W_1$, $w_2\in W_2$ and $w=w_1w_2$. Repeating literally the proof of \cite[Proposition 1.6]{EliseevIgnatyev}, we obtain the following: $$d_w=d_{w_1}^1d_{w_2}^2,\ c_w=c_{w_1}^1c_{w_2}^2.$$ Thus, to prove Theorem~\ref{mtheo:non_red}, it suffice to check it for irreducible root systems of type $D_n$, because $\Cp[\htt]$ is a unique factorization domain.

ii) Now, let $G\cong G_1\times G_2$, where $G_1$, $G_2$ are reductive subgroups of $G$, $T_i=T\cap G_i$ is a maximal torus in $G_i$, $i=1,2$, and the root system of $G_i$ with respect to $T_i$ is isomorphic to $\Phi_i$. Then $B_i=B\cap G_i$ is a Borel subgroup in $G_i$ containing $T_i$. Denote by $\Fo_i=G_i/B_i$ the corresponding flag variety. Then $\Fo=\Fo_1\times\Fo_2$ and $T_p\Fo=T_p\Fo_1\times T_p\Fo_2$ as algebraic varieties. If $w\in W$ and $w=w_1w_2$, $w_i\in W_i$, $i=1,2$, then $\rtc{w}\cong\rtc{w_1,G_1}\times\rtc{w_2,G_2}$ as affine varieties. Here $\rtc{w_i,G_i}$, $i=1,2$, denotes the tangent cone to the Schubert subvariety $X_{w_i}$ of the flag variety $\Fo_i$. Furthermore, note that $w$ is an involution if and only if $w_1$ and $w_2$ are involutions, too. This means that it suffice to prove that Theorem~\ref{mtheo:red} holds for all irreducible root system of type $D_n$.}

\sect{Non-reduced tangent cones}\label{sect:non_red}

\sst Throughout this Section, $\Phi$ denotes an irreducible root system of type $D_n$, $n\geq4$. In this Subsection, we briefly recall some facts about $\Phi$. Let $\epsi_1$, $\ldots$, $\epsi_n$ be the standard basis of the Euclidean space~$\Rp^n$. As usual, we identify the set $\Phi^+$ of positive roots with the following subset of $\Rp^n$:
\begin{equation*}
D_n^+=\{\epsi_i-\epsi_j,~\epsi_i+\epsi_j,~1\leq i<j\leq n\},\\
\end{equation*}
so $W$ can be considered as a subgroup of the orthogonal group $O(\Rp^n)$.

Let $S_{\pm n}$ denote the symmetric group on $2n$ letters $1,\ldots,n,-n,\ldots,-1$. The Weyl group $W$ is isomorphic to the \emph{even-signed hyper\-octa\-hedral group}, that is, the subgroup of $S_{\pm n}$ consisting of permutations ${w\in S_{\pm n}}$ such that $w(-i)=-w(i)$ for all $1\leq i\leq n$, and $\#\{i>0\mid w(i)<0\}$ is even. The isomorphism is given by
\begin{equation*}
\begin{split}
&s_{\epsi_i-\epsi_j}\mapsto(i,j)(-i,-j),\\
&s_{\epsi_i+\epsi_j}\mapsto(i,-j)(-i,j).\\
\end{split}
\end{equation*}
Here $s_{\alpha}$ is the reflection in the hyperplane orthogonal to a root $\alpha$. In the sequel, we will identify $W$ with the even-signed hyperoctahedral group.

\nota{i) Note that every $w\in W$ is completely determined by its restriction to the subset $\{1,\ldots,n\}$. This allows us to use the usual two-line notation: if $w(i)=w_i$ for $1\leq i\leq n$, then we will write $w=\begin{pmatrix}1&2&\ldots&n\\w_1&w_2&\ldots&w_n\end{pmatrix}$. For instance, if $\Phi=D_5$, then $$s_{\epsi_1+\epsi_5}s_{\epsi_2+\epsi_4}s_{\epsi_2-\epsi_4}=\begin{pmatrix}1&2&3&4&5\\-5&-2&3&-4&-1\end{pmatrix}.$$

ii) Note also that the set of simple roots has the following form: $\Delta=\{\alpha_1,\ldots,\alpha_n\}$, where $\alpha_1=\epsi_1-\epsi_2$, $\ldots$, $\alpha_{n-1}=\epsi_{n-1}-\epsi_n$, and $\alpha_n=\epsi_{n-1}+\epsi_n$.}

We say that $v$ is less or equal to $w$ with respect to the \emph{Bruhat order}, written $v\leq w$, if some reduced decomposition for $v$ is a subword of some reduced decomposition for $w$. It is well-known that this order plays the crucial role in many geometric aspects of theory of algebraic groups. For instance, the Bruhat order encodes the incidences among Schubert varieties, i.e., $X_v$ is contained in $X_w$ if and only if $v\leq w$. It turns out that $c_{w,v}$ is non-zero if and only if $v\leq w$ \cite[Corollary 3.2]{Kumar}. For example, $c_w=c_{w,\id}$ is non-zero for \emph{any} $w$, because $\id$ is the smallest element of $W$ with respect to the Bruhat order. Note that given $v,w\in W$, there exists $g_{w,v}\in S=\Cp[\htt]$ such that
\begin{equation}
c_{w,v}=g_{w,v}\cdot\prod_{\alpha>0,~s_{\alpha}v\leq w}\alpha^{-1},\label{formula:dyer}
\end{equation}
see \cite{Dyer} and \cite[Theorem 7.1.11]{BilleyLakshmibai}

There exists a nice combinatorial description of the Bruhat order on the even-signed hyper\-octa\-hedral group. Namely, given $w\in W$, denote by $X_w$ the $2n\times 2n$ matrix of the form
\begin{equation*}
(X_w)_{i,j}=\begin{cases}1,&\text{ if }w(j)=i,\\
0&\text{ otherwise}.
\end{cases}
\end{equation*}
The rows and the columns of this matrix are indicated by the numbers $1,\ldots,n,-n,\ldots,1$. It is called the $0$--1 matrix, permutation matrix or rook placement for $w$. Define the matrix $R_w$ by putting its $(i,j)$th element to be equal to the rank of the lower left $(n-i+1)\times j$ submatrix of $X_w$. In other words, $(R_w)_{i,j}$ is just the number or rooks located non-strictly to the South-West from $(i,j)$.

\exam{Let $n=4$, $w=\begin{pmatrix}1&2&3&4\\-2&4&1&-3\end{pmatrix}$. Here we draw the matrices $X_w$ and $R_w$ (rooks are marked by $\otimes$):
\begin{equation*}X_w=
\mymatrix{
\pho& \pho& \otimes& \pho& \pho& \pho& \pho& \pho\\
\pho& \pho& \pho& \pho& \pho& \pho& \pho& \otimes\\
\pho& \pho& \pho& \pho& \otimes& \pho& \pho& \pho\\
\pho& \otimes& \pho& \pho& \pho& \pho& \pho& \pho\\
\pho& \pho& \pho& \pho& \pho& \pho& \otimes& \pho\\
\pho& \pho& \pho& \otimes& \pho& \pho& \pho& \pho\\
\otimes& \pho& \pho& \pho& \pho& \pho& \pho& \pho\\
\pho& \pho& \pho& \pho& \pho& \otimes& \pho& \pho\\}
\ ,\ R_w=
\mymatrix{
\fho1& \fho2& \fho3& \fho4& \fho5& \fho6& \fho7& \fho8\\
\fho1& \fho2& \fho2& \fho3& \fho4& \fho5& \fho6& \fho7\\
\fho1& \fho2& \fho2& \fho3& \fho4& \fho5& \fho6& \fho6\\
\fho1& \fho2& \fho2& \fho3& \fho3& \fho4& \fho5& \fho5\\
\fho1& \fho1& \fho1& \fho2& \fho2& \fho3& \fho4& \fho4\\
\fho1& \fho1& \fho1& \fho2& \fho2& \fho3& \fho3& \fho3\\
\fho1& \fho1& \fho1& \fho1& \fho1& \fho2& \fho2& \fho2\\
\fho0& \fho0& \fho0& \fho0& \fho0& \fho1& \fho1& \fho1\\
}\ .
\end{equation*}\label{exam:Bruhat}}

Let $w\in W$. Given $a,b\in\{1,2,\ldots,n\}$, we say that $\empr{a}{b}$ is an \emph{empty rectangle} for $w$, if $$\{i\in[\pm n]\mid |i|\geq b\text{ and }|w(i)|\geq a\}=\varnothing.$$ Here $[\pm n]=\{1,\ldots,n,-n,\ldots,-1\}$. For instance, in the previous example $\empr{4}{3}$ and $\empr{4}{4}$ are empty rectangles for $w$. Let $X$ and $Y$ be matrices with integer entries. We say that $X\leq Y$ if $X_{i,j}\leq Y_{i,j}$ for all $i,j$. It turns out that given $v$, $w\in W$, $v\leq w$ if and only if
\begin{equation}
\begin{split}
&\text{i) }R_v\leq R_w;\\
&\text{ii) for all $a,b\in\{1,\ldots,n\}$, if $\empr{a}{b}$ is an empty rectangle for both $v$ and $w$}\\
&\hphantom{\text{ii) }}\text{and $(R_v)_{-(a-1),b-1}=(R_w)_{-(a-1),b-1}$, then $(R_v)_{-(a-1),n}\equiv(R_w)_{-(a-1),n}\pmod{2}$.}\\
\end{split}\label{formula:Bruhat}
\end{equation}(See, e.g., \cite[Theorem 8.2.8]{BjornerBrenti}.)

\sst\label{sst:support} In this Subsection, we introduce some more notation and prove technical, but crucial\break Lemma~\ref{lemm:crucial_non_red}. We define the maps $\row\colon\Phi^+\to\Zp$ and $\col\colon\Phi^+\to\Zp$ by
\begin{equation*}
\begin{split}
&\row(\epsi_i-\epsi_j)=j,~\row(\epsi_i+\epsi_j)=-j,\\
&\col(\epsi_i-\epsi_j)=\col(\epsi_i+\epsi_j)=i.
\end{split}
\end{equation*}
For any $k\in[\pm n]$, put
\begin{equation*}
\begin{split}
\Ro_k&=\{\alpha\in\Phi^+\mid\row(\alpha)=k\},\\
\Co_k&=\{\alpha\in\Phi^+\mid\col(\alpha)=k\}.\\
\end{split}
\end{equation*}
The set $\Ro_k$ (resp. $\Co_k$) is called the $k$th \emph{row} (resp. the $k$th \emph{column}) of $\Phi^+$.

\defi{An involution $w\in W$ is\label{defi:basic_involution} called \emph{basic}, if $$\{i\in\{1,\ldots,n\}\mid w(i)=-i\}=\varnothing.$$}

\defi{Let $\sigma\in W$ be a basic involution.\label{defi:support} We define the \emph{support} $\Supp{\sigma}$ of the involution~$\sigma$ by the following rule:
\begin{equation*}
\begin{split}
&\text{if }1\leq i<j\leq n\text{ and }\sigma(i)=j,\text{ then }\epsi_i-\epsi_j\in\Supp{\sigma},\\
&\text{if }1\leq i<j\leq n\text{ and }\sigma(i)=-j,\text{ then }\epsi_i+\epsi_j\in\Supp{\sigma}.\\
\end{split}
\end{equation*}
By definition, $\Supp{\sigma}$ is an orthogonal subset of $\Phi^+$. Note that $$\sigma=\prod_{\beta\in\Supp{\sigma}}s_{\beta},$$ where the product is taken in any fixed order. Note that for any $k$ one has $$|\Supp{\sigma}\cap\Co_k|\leq1,~|\Supp{\sigma}\cap\Ro_k|\leq1.$$ Note also that if $w$ is not basic, then, in general, there are several different ways to define $\Supp{w}$, see Remark~\ref{nota:why_basic_red}~(ii) below.}

\exam{Let $\Phi=D_6$ and $\sigma=\begin{pmatrix}1&2&3&4&5&6\\-6&2&5&4&3&-1\end{pmatrix}$. Then $$\Supp{\sigma}=\{\epsi_1+\epsi_6,\epsi_3-\epsi_5\}.$$}

\nota{i) Denote the set of involutions (resp. of basic involutions) by $I(W)$ (resp. by $B(W)$). By \cite[Proposition 2.3]{Ignatyev2}, if $\sigma,\tau\in I(W)$, then
\begin{equation}
R_{\sigma}\leq R_{\tau}\text{ if and only if }R_{\sigma}^*\leq R_{\tau}^*,\label{formula:Bruhat_invol}
\end{equation}
where $R_w^*$ is the strictly lower-triangular part of $R_w$, i.e.,
\begin{equation*}
(R_w^*)_{i,j}=\begin{cases}
(R_w)_{i,j}&\text{if }i>j,\\
0,&\text{if }i\leq j.
\end{cases}
\end{equation*}

ii) Using Formulas (\ref{formula:Bruhat}) or (\ref{formula:Bruhat_invol}), one can easily check that if $\alpha\in\Co_1$ and $\beta\notin\Co_1$, then $s_{\alpha}\nleq s_{\beta}$. One can also check that
\begin{equation*}
s_{\epsi_1-\epsi_2}<\ldots<s_{\epsi_1-\epsi_n},~s_{\epsi_1+\epsi_n}<\ldots<s_{\epsi_1+\epsi_2}.
\end{equation*}
Further, $s_{\epsi_1-\epsi_i}<s_{\epsi_1+\epsi_j}$ for all $i,j\in\{1,\ldots,n\}$ such that $i<n$ or $j<n$, but $s_{\epsi_1-\epsi_n}\nless s_{\epsi_1+\epsi_n}$ and $s_{\epsi_1+\epsi_n}\nless s_{\epsi_1-\epsi_n}$.}

The following Lemma plays the crucial role in the proof of Theorem~\ref{mtheo:non_red} (cf. \cite[Lemmas 2.4, 2.5]{EliseevIgnatyev} and \cite[Lemma 2.6]{BochkarevIgnatyevShevchenko}).

\lemmp{Let $w\in W$ be a basic involution. If $\Supp{w}\cap\Co_1=\varnothing$\textup{,} then $\alpha$ divides $d_w$ in the polynomial ring $\Cp[\htt]$ for all $\alpha\in\Co_1$. If $\Supp{w}\cap\Co_1=\{\beta\}$\textup{,} then $\beta$ does not divide $d_w$ in $\Cp[\htt]$.\label{lemm:crucial_non_red}}{Denote by $\wt W$ the subgroup of $W$ generated by $s_2,\ldots,s_n$. Suppose $\Supp{w}\cap\Co_1=\varnothing$, then $w\in\wt W$. Denote by $\wt\Phi$ the root system corresponding to $\wt W$; in fact, $\wt\Phi^+=\Phi^+\setminus\Co_1$.

Let $\wt d_w\in\wt S=\Cp[\alpha_2,\ldots,\alpha_n]$ be the Kostant--Kumar polynomial of $w$ considered as an element of~$\wt W$; define $\wt c_w\in\Cp(\alpha_2,\ldots,\alpha_n)$ by the similar way. Since $\wt W$ is a parabolic subgroup of $W$, the length of $w$ as an element of $\wt W$ equals the length of $w$ as an element of $W$. Further, any reduced decomposition for $w$ in $\wt W$ is a reduced decomposition for $w$ in $W$. This means that $\wt c_w=c_w$, so
\begin{equation*}
\begin{split}
d_w=(-1)^{l(w)}\cdot\prod_{\alpha\in\Phi^+}\alpha\cdot c_w=(-1)^{l(w)}\cdot\prod_{\alpha\in\Co_1}\alpha\cdot\prod_{\alpha\in\wt\Phi^+}\alpha\cdot\wt c_w=\wt d_w\cdot\prod_{\alpha\in\Co_1}\alpha.
\end{split}
\end{equation*}
In particular, $\alpha$ divides $d_w$ for all $\alpha\in\Co_1$.

Now, suppose $\Supp{w}\cap\Co_1=\{\beta\}$. By \cite[Proposition 1.10]{Humpreys2}, there exists a unique $v\in\wt W$ such that $w=uv$ and $l(us_i)=l(u)+1$ for all $2\leq i\leq n$ (or, equivalently, $u(\alpha_i)>0$ for all $2\leq i\leq n$). Furthermore, $l(w)=l(u)+l(v)$. One can easily check that
\begin{equation*}\predisplaypenalty=0
\begin{split}
&\text{if }\beta=\epsi_1-\epsi_j\text{ (i.e., $w(1)=j$), then}\\
&u=s_{j-1}\ldots s_2s_1\\
&\hphantom{u}=\begin{cases}\begin{pmatrix}1&2&3&\ldots&j-1&j&j+1&\ldots&n-1&n\\j&1&2&\ldots&j-2&j-1&j+1&\ldots&n-1&n\\\end{pmatrix},&\text{if }j<n,\\
\begin{pmatrix}1&2&3&\ldots&n-1&n\\n&1&2&\ldots&n-2&n-1\\\end{pmatrix},&\text{if }j=n,\\
\end{cases}\\
&\text{if }\beta=\epsi_1+\epsi_j\text{ (i.e., $w(1)=-j$), then}\\
&u=s_js_{j+1}\ldots s_{n-1}s_n s_{n-2}s_{n-3}\ldots s_2s_1\\
&\hphantom{u}=\begin{cases}
\begin{pmatrix}1&2&3&\ldots&j-1&j&j+1&\ldots&n-1&n\\-j&1&2&\ldots&j-2&j-1&j+1&\ldots&n-1&-n\\\end{pmatrix},&\text{if }j<n,\\
\begin{pmatrix}1&2&3&\ldots&n-1&n\\-n&1&2&\ldots&n-2&-(n-1)\\\end{pmatrix},&\text{if }j=n.\\
\end{cases}
\end{split}
\end{equation*}

For instance, consider the case $\beta=\epsi_1+\epsi_j$ (the case $\beta=\epsi_1-\epsi_j$ can be considered similarly). Recall that $W$ acts on $\Cp(\htt)$ by automorphisms. Using (\ref{formula_x_v_x_w}) and arguing as in the proof of
\cite[Lemma 2.5]{EliseevIgnatyev}, one can easily show that
\begin{equation}
c_w=-\dfrac{c_{us_1,g_0}g_0(c_{v,g_0^{-1}})}{\beta}-\sum_{g\leq u,\ g^{-1}\leq v,\ g\neq g_0}\dfrac{c_{us_1,g}g(c_{v,g^{-1}})}{g(\alpha_1)}=
\beta^{-1}\cdot g_0(c_{v,g_0^{-1}})\cdot\dfrac{K}{L}+\dfrac{M}{N}\label{formula:KLMN}
\end{equation}
(cf. Formula (7) from \cite{EliseevIgnatyev}). Here
\begin{equation*}
\begin{split}
g_0&=us_1=s_js_{j+1}\ldots s_{n-1}s_n s_{n-2}s_{n-3}\ldots s_2\\
&=\begin{cases}\begin{pmatrix}1&2&3&\ldots&j-1&j&j+1&\ldots&n-1&n\\1&-j&2&\ldots&j-2&j-1&j+1&\ldots&n-1&-n\\\end{pmatrix},&\text{if }j<n,\\
\begin{pmatrix}1&2&3&\ldots&n-1&n\\1&-n&2&\ldots&n-2&-(n-1)\\\end{pmatrix},&\text{if }j=n,\\
\end{cases}
\end{split}
\end{equation*}
and $K,L$ and $M,N\in\Cp[\htt]$ are pairs of coprime polynomials such that $\beta$ divides neither $K$ nor $N$.

To prove that $\beta$ does not divide $d_w$, it is enough to show that $c_{v,g_0^{-1}}\neq0$, i.e., $v\geq g_0^{-1}$ (or, equivalently, $v^{-1}\geq g_0$). Arguing as in the proof of \cite[Lemma 2.6]{BochkarevIgnatyevShevchenko}, we obtain $R_{v^{-1}}\geq R_{g_0}$. Thus, it remains to check that the second condition in the definition of the Bruhat order is satisfied. Suppose that $\empr{a}{b}$ is an empty rectangle for both $v^{-1}$ and $g_0$ and $(R_{g_0})_{-(a-1),b-1}=(R_{v^{-1}})_{-(a-1),b-1}$. We must prove that $(R_{g_0})_{-(a-1),n}\equiv(R_{v^{-1}})_{-(a-1),n}\pmod{2}$.

If $j<n$, then $g_0(n)=-n$, so there are no empty rectangles for~$g_0$, hence $j=n$. In this case, $a=n$ and $b\geq 3$. For example, on the picture below we draw $X_{g_0}$ for $n=5$, $b=4$. Entries from the empty rectangle $\empr{5}{4}$ are grey.
\begin{equation*}\predisplaypenalty=0
\mymatrix{
\otimes& \pho& \pho& \pho& \pho& \pho& \pho& \pho& \pho& \pho\\
\pho& \pho& \otimes& \pho& \pho& \pho& \pho& \pho& \pho& \pho\\
\pho& \pho& \pho& \otimes& \pho& \pho& \pho& \pho& \pho& \pho\\
\pho& \pho& \pho& \pho& \pho& \otimes& \pho& \pho& \pho& \pho\\
\pho& \pho& \pho& \gray\pho& \gray\pho& \gray\pho& \gray\pho& \pho& \otimes& \pho\\
\pho& \otimes& \pho& \gray\pho& \gray\pho& \gray\pho& \gray\pho& \pho& \pho& \pho\\
\pho& \pho& \pho& \pho& \otimes& \pho& \pho& \pho& \pho& \pho\\
\pho& \pho& \pho& \pho& \pho& \pho& \otimes& \pho& \pho& \pho\\
\pho& \pho& \pho& \pho& \pho& \pho& \pho& \otimes& \pho& \pho\\
\pho& \pho& \pho& \pho& \pho& \pho& \pho& \pho& \pho& \otimes\\}
\end{equation*}
Clearly, $(R_{g_0})_{-(n-1),b-1}=0$, hence $(R_{v^{-1}})_{-(n-1),b-1}=0$. At the same time, $(R_{g_0})_{-(n-1),n}=1$, so we must check that $(R_{v^{-1}})_{-(n-1),n}$ is odd.

By definition,
$$(R_{v^{-1}})_{-(n-1),n}=\#\{i\in\{1,\ldots,n\}\mid v^{-1}(i)\in\{-1,\ldots,-(n-1)\}\}.$$
On the other hand, $v^{-1}(2)=wu(2)=w(1)=-n$, hence $\#\{i\in\{1,\ldots,n\}\mid v^{-1}(i)=-n\}=\#\{n\}=1$. Since the number
\begin{equation*}
(R_{v^{-1}})_{-(n-1),n}+1=\#\{i\in\{1,\ldots,n\}\mid v^{-1}(i)<0\}
\end{equation*}
is even (by definition of $W$), we conclude that $(R_{v^{-1}})_{-(n-1),n}$ is odd, as required.}

\sst Things now are ready for the proof of our first main result, Theorem~\ref{mtheo:non_red}. The proof immediately follows from the Proposition~\ref{prop:non_red_eq_C_1} below (cf. \cite[Propositions 2.6, 2.7, 2.8]{EliseevIgnatyev} and \cite[Propositions 2.7, 2.8]{BochkarevIgnatyevShevchenko}). Our goal is to check that if $\sigma,\tau$ are distinct basic involutions in $W$, then their Kostant--Kumar polynomials do not coincide, and, consequently, the tangent cones $C_{\sigma}$ and $C_{\tau}$ do not coincide as subschemes of $T_p\Fo$. We will proceed by induction on $n$ (the base is trivial).

\propp{Let $\sigma,\tau\in W$ be distinct basic involutions\label{prop:non_red_eq_C_1}. Then $d_{\sigma}\neq d_{\tau}$.}{
If $\Supp{\sigma}\cap\Co_1\neq\Supp{\tau}\cap\Co_1$, then one can repeat literally the proof of \cite[Pro\-po\-si\-tion~2.7]{BochkarevIgnatyevShevchenko} to obtain the result. Namely, if $\Supp{\sigma}\cap\Co_1=\{\beta\}$ and $\Supp{\tau}\cap\Co_1=\varnothing$, then $\beta$ does not divide $d_{\sigma}$ by the previous Lemma. But, thanks to formula (\ref{formula:dyer}), $\beta$ divides $d_{\tau}$, so $d_{\sigma}\neq d_{\tau}$. On the other hand, suppose that $\Supp{\sigma}\cap\Co_1={\beta}$, $\Supp{\tau}\cap\Co_1=\beta'$, $\beta\nless\beta'$, then $\beta$ divides $d_{\tau}$ (by formula (\ref{formula:dyer})), but $\beta$ does not divide $d_{\sigma}$ (by the previous Lemma), so $d_{\sigma}\neq d_{\tau}$.

From now on, we may assume that $\Supp{\sigma}\cap\Co_1=\Supp{\tau}\cap\Co_1$. If $\Supp{\sigma}\cap\Co_1=\Supp{\tau}\cap\Co_1=\varnothing$, then the inductive assumption completes the proof. Suppose $\Supp{\sigma}\cap\Co_1=\Supp{\tau}\cap\Co_1=\{\beta\}$. Let $u$ be as in the proof of Lemma~\ref{lemm:crucial_non_red}. There are two cases:
\begin{equation*}
\begin{split}
&\text{i) }\beta=\epsi_1-\epsi_j,\text{ i.e., }w(1)=j,\\
&\text{ii) }\beta=\epsi_1+\epsi_j,\text{ i.e., }w(1)=-j.\\
\end{split}
\end{equation*}

If $\beta=\epsi_1-\epsi_j$, then one can repeat literally the proof of Case (i) of \cite[Proposition 2.8]{BochkarevIgnatyevShevchenko}, so we may assume that $\beta=\epsi_1+\epsi_j$. Arguing as in the proof of Case (ii) of \cite[Proposition 2.8]{BochkarevIgnatyevShevchenko}, we obtain that $$c_{v,g_0^{-1}}=c_{v_2,\id}\cdot\prod_{i=3}^{n-1}(\epsi_2-\epsi_i)^{-1}\cdot
\prod_{i=j+1}^{n-1}(\epsi_2+\epsi_i)^{-1}\cdot(\epsi_2+\epsi_n)^{-2}.$$
Here $w=aw_2a^{-1}$, $a=s_2s_3\ldots s_{n-2}s_ns_{n-1}\ldots s_{j+1}s_j$, $w_2=u_2v_2$, $\Supp{w_2}\cap\Co_1=\{\alpha_1\}$, $u_2=s_1$, and $v_2\in\wt W$ is an involution.

Now, consider the involutions $\sigma$ and $\tau$. Put $\sigma=uv_{\sigma}$, $\tau=uv_{\tau}$, where $u$ is as above. Put also $\sigma=a\sigma_2a^{-1}$, $\tau=a\tau_2a^{-1}$, $\sigma_2=u_2v_{\sigma}^2$, $\tau_2=u_2v_{\tau}^2$, where $u_2=s_1$. By the inductive assumption, $c_{v_{\sigma}^2,\id}\neq c_{v_{\tau}^2,\id}$, hence $c_{v_{\sigma},g_0^{-1}}\neq c_{v_{\tau},g_0^{-1}}$. Arguing as in the last two paragraphs of the proof of \cite[Proposition 2.8]{EliseevIgnatyev}, one can conclude the proof.

Namely, one can easily deduce from formula (\ref{formula:KLMN}) that if $c_\sigma=c_{\tau}$, then $\beta$ divides $P_{\sigma}Q_{\tau}-P_{\tau}Q_{\sigma}$, where $P_{\sigma}$ and $Q_{\sigma}$ (resp. $P_{\tau}$ and $Q_{\tau}$) are coprime polynomials such that $g_0(c_{v_{\sigma},g_0^{-1}})=P_{\sigma}/Q_{\sigma}$ (resp. $g_0(c_{v_{\tau},g_0^{-1}})=P_{\tau}/Q_{\tau}$). But these polynomials belong to the subalgebra of $\Cp[\htt]$ generated by $\alpha_2,\ldots,\alpha_n$, so $c_{v_{\sigma},g_0^{-1}}=c_{v_{\tau},g_0^{-1}}$, a contradiction.}

\sect{Reduced tangent cones}\label{sect:red}

\sst In this\label{sst:red_coadjoint} Section we will prove our second main result, Theorem~\ref{mtheo:red}. Throughout the Section, we will assume that every $\Phi$ is of type $D_n$, $n\geq4$. In this Subsection, we briefly describe connections between tangent cones and coadjoint orbits of $U$, the unipotent radical of the Borel subgroup $B$.

Denote by $\gt$, $\bt$, $\nt$ the Lie algebras of $G$, $B$, $U$ respectively, then $T_p\Fo$ is naturally isomorphic to the quotient space $\gt/\bt$. Using the Killing form on $\gt$, one can identify the latter space with the dual space~$\nt^*$. The group $B$ acts on $\Fo$ by conjugation. Since $p$ is $B$-stable, $B$ acts on the tangent space $T_p\Fo\cong\nt^*$. This action is called \emph{coadjoint}. We denote the result of coadjoint action by $b.\lambda$, $b\in B$, $\lambda\in\nt^*$. In 1962, A.A. Kirillov discovered that orbits of this action play an important role in representation theory of $B$ and $U$, see, e.g., \cite{Kirillov1}, \cite{Kirillov2}. We fix a basis $\{e_{\alpha},~\alpha\in\Phi^+\}$ of $\nt$ consisting of root vectors. Let $\{e_{\alpha}^*,~\alpha\in\Phi^+\}$ be the dual basis of~$\nt^*$. Let $w\in W$ be a basic involution. Put $$f_w=\sum_{\beta\in\Supp{w}}e_{\beta}^*\in\nt^*.$$

\defi{We say that the $B$-orbit $\Omega_w$ and the $U$-orbit $\Theta_w$ of $f_w$ are \emph{associated} with the involution $w$.}

One can easily check that $\Theta_w\subset\Omega_w\subseteq\rtc{w}$. Further, $\rtc{w}$ is $B$-stable (in fact, the tangent cone to an arbitrary Schubert variety is $B$-stable). Orbits associated with involutions were studied by A.N.~{Pa\-nov}~\cite{Panov} and the second author \cite{Ignatyev1}, \cite{Ignatyev2}, \cite{Ignatyev3}, \cite{Ignatyev4} (see also the Kostant's papers \cite{Kostant1}, \cite{Kostant2}, \cite{Kostant2} for the connections with the center of enveloping algebra of $\nt$). In particular, it was shown in \cite[Theorem 1.2]{Ignatyev3} that
\begin{equation}
\dim\Theta_w=l(w)-|\Supp{w}|.\label{formula:dim_Theta}
\end{equation}
We need the following corollary of this fact (cf. \cite[Pro\-po\-si\-tion~4.1]{Ignatyev1} and \cite[Theorem 3.1]{Ignatyev2}).
\lemmp{If $w\in W$ is a basic involution, then
\begin{equation}
\dim\Omega_w=l(w).\label{formula:dim_Omega}
\end{equation}}{Denote $D=\Supp{w}$. Let $\xi\colon D\to\Cp^{\times}$ be a map. Denote by $\Theta_{w,\xi}$ the $U$-orbit of the linear form $$f_{w,\xi}=\sum_{\beta\in D}\xi(\beta)e_{\beta}^*.$$
In particular, $f_w=f_{w,\xi_0}$, where $\xi_0(\beta)=1$ for all $\beta\in D$.

Without loss of generality, we can identify $G$ with the group $\SO_{2n}(\Cp)$ of all invertible $2n\times 2n$ matrices $g$ of determinant 1 such that $g^tJg=J$, where $J$ is the symmetric $2n\times 2n$ with 1's on the antidiagonal and 0's elsewhere. Then $T$ (resp. $B$ and $U$) is the group of all diagonal (resp. upper-triangular and upper-triangular with 1's on the diagonal) matrices from $G$. Moreover, $\gt$ is the algebra of $2n\times2n$ matrices $x$ of zero trace satisfying $x^tJ+Jx=0$, and $\htt$ (resp. $\bt$ and $\nt$) is the algebra of all diagonal (resp. upper-triangular and upper-triangular with 0's on the diagonal) matrices from $\gt$. Using Killing form of $\gt$, one can identify $\nt^*$ with the space $\nt^t$ of all lower-triangular matrices from $\gt$ with 0's on the diagonal. Under this identification, the coadjoint action of $B$ has a simple form
\begin{equation}b.\lambda=(b\lambda b^{-1})_{\low},~b\in B,~\lambda\in\nt^*,\label{formula:coad_explicit}
\end{equation}
where $A_{\low}$ denotes the strictly lower-triangular part of a matrix $A$.

First, we claim that if $\xi_1\neq \xi_2$, then $\Theta_{w,\xi_1}\neq\Theta_{w,\xi_2}$. Indeed, let $\wt U$ be the group of all $2n\times 2n$ upper-triangular matrices with 1's on the diagonal. This group acts on the space $\wt\nt$ of all upper-triangular $2n\times2n$ matrices with 0's on the diagonal by the adjoint action, hence one can consider the dual (coadjoint) action of this group on the space $\wt\nt^*$. Using Killing form of $\mathfrak{gl}_{2n}(\Cp)$, one can identify $\wt\nt^*$ with the space $\wt\nt^t$ of all lower-triangular $2n\times2n$ matrices with 0's on the diagonal. Under this identification, the coadjoint action of~$\wt U$ is given again by formula (\ref{formula:coad_explicit}). Let $\wt\Theta_{w,\xi}\subset\wt\nt^*$ be the $\wt U$-orbit of $f_{w,\xi}$, then, clearly, $\Theta_{w,\xi}\subseteq\wt\Theta_{w,\xi}$ for any $\xi$. Since $w$ is an involution in $S_{\pm n}$, it follows from \cite[Theorem 1.4]{Panov} that $\wt\Theta_{w,\xi_1}\neq\wt\Theta_{w,\xi_2}$. Thus, $\Theta_{w,\xi_1}\neq\Theta_{w,\xi_2}$, as required.

Second, we claim that $\Omega_w=\bigcup_{\xi}\Theta_{w,\xi}$, where the union is taken over all maps from $D$ to $\Cp^{\times}$. Indeed, is is well-known that the \emph{exponential map} $$\exp\colon\nt\to U,~x\mapsto\sum_{i=0}^{\infty}\dfrac{x^i}{i!}$$ is an isomorphism of affine varieties. Given $\alpha\in\Phi^+$, $s\in\Cp^{\times}$, put
\begin{equation*}
\begin{split}
&x_{\alpha}(s)=\exp{se_{\alpha}}=1+se_{\alpha},~x_{-\alpha}(s)=x_{\alpha}(s)^t,\\
&w_{\alpha}(s)=x_{\alpha}(s)x_{-\alpha}(-s^{-1})x_{\alpha}(s),~h_{\alpha}(s)=w_{\alpha}(s)w_{\alpha}(1)^{-1}.
\end{split}
\end{equation*}
Note that $h_{\alpha}(s)$ belongs to $T$.

Let $\xi\colon D\to\Cp^{\times}$ be a map, $\alpha\in D$ be a root. To any number $s\in\Cp^{\times}$, denote by $\sqrt{s}$ a complex number such that $\left(\sqrt{s}\right)^2=s$. One can trivially check by direct matrix calculations that $$h_{\alpha}(\sqrt{s}).f_{w,\xi}=s\xi(\alpha)e_{\alpha^*}+\sum_{\beta\in D,~\beta\neq\alpha}\xi(\beta)e_{\beta}^*.$$
Thus, $$\left(\prod_{\alpha\in D}h_{\alpha}\left(\sqrt{\xi(\alpha)}\right)\right).f_w=f_{w,\xi},$$
so $\Theta_{w,\xi}\subset\Omega_w$.

On the other hand, $B=U\rtimes T$ as algebraic groups. Since $T$ is generated by $h_{\alpha}(s)$, $\alpha\in\Phi^+$, $s\in\Cp^{\times}$, we see that if $h\in T$, then $h.f_{w,\xi}=f_{w,\xi'}$ for some map $\xi'\colon D\to\Cp^{\times}$. Thus, if $g\in B$ and $g=uh$, $u\in U$, $h\in T$, then $g.f_w=u.f_{w,\xi}$ for some $\xi$, so $\Omega_w=\bigcup_{\xi}\Theta_{w,\xi}$, as required.

Third, let $Z_B$ (resp. $Z_U$ and $Z_T$) be the stabilizer of $f_w$ under the coadjoint action of $B$ (resp. of~$U$ and $T$). Then
\begin{equation*}
\begin{split}
&\dim\Omega_w=\dim B-\dim Z_B,\\
&\dim\Theta_w=\dim U-\dim Z_U.
\end{split}
\end{equation*}
If $g=uh\in Z_B$, $u\in U$, $h\in T$, then $$g.f_w=u.(h.f_w)=u.f_{w,\xi}$$ for some $\xi$. If $f_w\neq f_{w,\xi}$, then $\Theta_w\neq\Theta_{w,\xi}$. Hence $f_w=f_{w,\xi}$, so $h\in Z_T$ and $u\in Z_U$. It follows that the map $$Z_U\times Z_T\to Z_B\colon(u,h)\mapsto uh$$
is an isomorphism of algebraic varieties, so $$\dim Z_B=\dim Z_U+\dim Z_T.$$

Finally, it follows that $X=\bigcup_{\xi}\{f_{w,\xi}\}$ is the $T$-orbit of $f_w$ (the union is taken over all maps from~$D$ to $\Cp^{\times}$). Thus, using (\ref{formula:dim_Theta}), we conclude that
\begin{equation*}
\begin{split}
\dim\Omega_w&=\dim B-\dim Z_B\\
&=\dim U+\dim T-\dim Z_U-\dim Z_T\\
&=\dim\Theta_w+\dim X=l(w)-|D|+|D|=l(w).
\end{split}
\end{equation*}
The proof is complete.}

\nota{i) Since $\dim\rtc{w}=\dim X_w=l(w)$, we conclude that $\overline{\Omega}_w$, the closure of $\Omega_w$, is an irreducible component of $\rtc{w}$ of maximal dimension. (In fact, $\rtc{w}$ is equidimensional.)

ii) If $w$ is not basic,\label{nota:why_basic_red} then there are several different ways to define $\Supp{w}$. For example, if $n=4$ and $w=\begin{pmatrix}1&2&3&4\\-1&-2&-3&-4\end{pmatrix}$, then there are three subsets $D\subset\Phi^+$ such that $w=\prod_{\beta\in D}s_{\beta}$:
\begin{equation*}
\begin{split}
&\{\epsi_1-\epsi_2,\epsi_1+\epsi_2,\epsi_3-\epsi_4,\epsi_3+\epsi_4\},\\
&\{\epsi_1-\epsi_3,\epsi_1+\epsi_3,\epsi_2-\epsi_4,\epsi_2+\epsi_4\},\\
&\{\epsi_1-\epsi_4,\epsi_1+\epsi_4,\epsi_2-\epsi_3,\epsi_2+\epsi_3\}.\\
\end{split}
\end{equation*}
(For some reasons, the first candidate is ``the best'', see \cite{Deodhar}, \cite{Springer}.)

So, one can define $\Theta_w$ and $\Omega_w$ using one of the definitions of the support of $w$. But there is no chance that formula (\ref{formula:dim_Omega}) holds for all non-basic involutions. Indeed, one can repeat literally the proof of the previous Lemma to obtain $\dim\Omega_w=\dim\Theta_w+|\Supp{w}|$. But if $w$ is not basic, then the dimension of $\dim\Theta_w$ can be strictly less than $l(w)-|\Supp{w}|$, see \cite[Theorem 1.2]{Ignatyev3}. That's why we restrict our attention to the case of basic involutions.}

Now, assume that $G'$ is a reductive subgroup of $G''$, $T'$ (resp. $T''$) is a maximal torus of $G'$ (resp. of $G''$), $T'=T''\cap G'$, $B'$ (resp.~$B''$) is a Borel subgroup of $G'$ (resp. of $G''$) containing $T'$ (resp.~$T''$), $B'=B''\cap G'$, and $\Phi'$ (resp. $\Phi''$) is the root system of $G'$ (resp. of $G''$) with respect to $T'$ (resp. to~$T''$). We denote by $W'$ (resp. by $W''$) the Weyl group of $\Phi'$ (resp. of $\Phi''$). Denote by $\Fo'=G'/B'$, $\Fo''=G''/B''$ the flag varieties. Put $p'=eB'\in\Fo'$, $p''=eB''\in\Fo''$. Let $U'$ (resp. $U''$) be the unipotent radical of $B'$ (resp. of $B''$), $U'=U''\cap B'$. Denote also by $\gt'$,~$\bt'$,~$\nt'$ the Lie algebras of $G'$, $B'$, $U'$ respectively. Define $\gt''$, $\bt''$, $\nt''$ by the similar way. One can consider the dual space $\nt'^*\cong\gt'/\bt'$ as a~subspace of $\nt''^*\cong\gt''/\bt''$. Hence we can consider $T_{p'}\Fo'$ as a subspace of $T_{p''}\Fo''$.

Pick involutions $w_1,w_2\in W'$. Let $C_i'$ be the reduced tangent cone at the point $p'$ to the Schubert subvariety $X_{w_i}'$ of the flag variety $\Fo'$, $i=1,2$. Similarly, let $C_i''$ be the reduced tangent cone at $p''$ to the Schubert subvariety $X_{w_i}''$ of $\Fo''$, $i=1,2$. Denote by $l'$ (resp. by $l''$) the length function on the Weyl group $W'$ (resp. on $W''$). Assume $C_1'=C_2'$. This implies that $$l'(w_1)=l'(w_2).$$ Note that $C_i'\subseteq C_i''$, hence $B''.C_i'\subseteq C_i''$, $i=1,2$. Denote by $\Omega_{w_i}'\subseteq\nt'^*$ the coadjoint $B'$-orbit associated with the involution $w_i$, $i=1,2$; define $\Omega_{w_i}''$ by the similar way. It follows from formula (\ref{formula:dim_Omega}) that
\begin{equation*}
\begin{split}
l''(w_i)&=\dim C_i''\geq\dim B''.C_i'\geq\dim B''.\Omega_{w_i}'\\
&=\dim\Omega_{w_i}''=l''(w_i),
\end{split}
\end{equation*}
because $\Omega_{w_i}''=B''.\Omega_{w_i}'$. This implies $l''(w_i)=\dim C_i''=\dim B''.C_i'$. But $C_1'=C_2'$, thus $\dim C_1''=\dim C_2''$. We obtain the following result:
\begin{equation}
\text{if $C_1'=C_2'$, then $l''(w_1)=l''(w_2)$.}\label{formula:if_cones_then_ls}
\end{equation}

\sst In this Subsection, \label{sst:red_D_n}we prove Theorem~\ref{mtheo:red}: if $w_1$, $w_2$ are basic involutions in the Weyl group~$W$ of type $D_n$, $n\geq4$, and $w_1\neq w_2$, then $\rtc{w_1}\neq\rtc{w_2}$ as subvarieties of $T_p\Fo$. Let $W''$ be of type $D_{n+2}$. Let
$$D_{n+2}^+=\{\eta_i-\eta_j,~\eta_i+\eta_j,~1\leq i<j\leq n+2\},$$
where $\{\eta_i\}_{i=1}^{n+2}$ is the standard basis of $\Rp^{n+2}$. Pick numbers $k_1,k_2$ such that $1\leq k_1<k_2\leq n+2$. Put $P=\{k_1,k_2\}$, $Q=\{1,\ldots,n+2\}\setminus P$, and
\begin{equation*}
\begin{split}
\wt W&=\{w\in W''\mid w(i)=i\text{ for all }i\in P\},\\
\wt W_2&=\{w\in W''\mid w(i)=i\text{ for all }i\in Q\}, \\
W'&=\{w\in W''\mid w(P)=P,~w(Q)=Q\}=\wt W\times\wt W_2.
\end{split}
\end{equation*}
Let $\Phi'$ (resp. $\wt\Phi$) be the root system of $W'$ (resp. of $\wt W$). Clearly, $\Phi'$ (resp. $\wt\Phi$) is of type $D_n\times A_1\times A_1$ (resp. of type $D_n$). Put $G''=\SO_{2n+4}(\Cp)$ and denote by $G'$ (resp. by $\wt G$) the subgroup of $G$ corresponding to~$\Phi'$ (resp. to $\wt\Phi$), then $G'\cong\SO_n(\Cp)\times\SO_2(\Cp)$. Put also
\begin{equation*}
\begin{split}
&A=\{1,\ldots,k_1-1\},\\
&B=\{k_1+1,\ldots,k_2-1\},\\
&C=\{k_2+1,\ldots,n+2\}.
\end{split}
\end{equation*}

Now, let $\Phi=D_n$. We can assume without loss of generality that $G=\SO_n(\Cp)$. We identify $\Phi$ with $\wt\Phi$ by the map $\epsi_k\mapsto\eta_{k'}$, where
\begin{equation*}
k'=\begin{cases}k,&\text{if }k\leq k_1-1,\\
k+1,&\text{if }k_1\leq k\leq k_2-2,\\
k+2,&\text{if }k_2-1\leq k\leq n.
\end{cases}
\end{equation*}
This identifies $G$ (resp. $W$) with $\wt G$ (resp. with $\wt W$). We denote the image in $\wt W$ of an element $w\in W$ under this identification by $\wt w$. Let $w\in W$ be an involution. Arguing as in the proof of \cite[Lemma~3.2]{BochkarevIgnatyevShevchenko}, we obtain the following result.

\mlemm{\textup{i)} If $w'=\wt ws_{\eta_{k_1}-\eta_{k_2}}$\textup{,} then the length of $w'$ in the \label{lemm:l_w_dash_A_n}Weyl group $W''$ equals $$l''(w')=2(k_2-k_1-1)+4|\wt w(A)\cap B^-|+4|\wt w(A)\cap A^-|+4|\wt w(A)\cap C^{\pm}|+l(w)+1.$$
\textup{ii)} If $w'=\wt ws_{\eta_{k_1}+\eta_{k_2}}$\textup{,} then
$$l''(w')=2(k_2-k_1-1)+4|\wt w(A)\cap A^-|+4|\wt w(A)\cap B^-|+4|C|+l(w)+1.$$}
(By a slight abuse of notation, here we consider $\wt w$ as an element of $S_{\pm(n+2)}$ and, at the same time, as an element of $\wt W$, i.e., as an element of $W''$ such that $\wt w(k_1)=k_1$ and $\wt w(k_2)=k_2$.)

\textsc{Proof of Theorem~\ref{mtheo:red}.} Assume $\rtc{w_1}=\rtc{w_2}$. In particular, $$l(w_1)=\dim\rtc{w_1}=\dim\rtc{w_2}=l(w_2).$$ Since $w_1\neq w_2$, there exists $1\leq k\leq n$ such that $w_1(\epsi_i)=w_2(\epsi_i)$ for $1\leq i\leq k-1$, and $w_1(\epsi_k)\neq w_2(\epsi_k).$ Assume without loss of generality that $w_1(\epsi_k)<w_2(\epsi_k)$, i.e., $w_2(\epsi_k)-w_1(\epsi_k)$ is a sum of positive roots. Note that $w_1(\epsi_k)\neq\pm\epsi_k$, because $w_1(\epsi_i)=w_2(\epsi_i)$ for all $i$ from 1 to $k-1$. Put $k_1=k+1$, so $A=\{1,\ldots,k\}$ and $\wt w_1(a)=\wt w_2(a)$ for all $a\in A\setminus\{k\}$. We consider three different cases.

i) Suppose $w_1(\epsi_k)<0$, $w_2(\epsi_k)>0$. Here we put $k_2=n+2$, so $C=\varnothing$ and $$(\wt w_i(A)\cap A^-)\cup(\wt w_i(A)\cap B^-)=\wt w_i(A)\cap\{-1,\ldots,-(n+2)\},~i=1,2.$$
Let $w_i'=\wt w_is_{\eta_{k_1}-\eta_{k_2}}$, $i=1,2$. Since $$\wt w_1(A)\cap\{-1,\ldots,-(n+2)\}=\wt w_2(A)\cap\{-1,\ldots,-(n+2)\}\cup\{k\},$$
Lemma~\ref{lemm:l_w_dash_A_n}~(i) shows that $l''(w_1')\neq l''(w_2')$. On the other hand, $\rtc{w_1}=\rtc{w_2}$ implies $C_1'=C_2'$, which contradicts (\ref{formula:if_cones_then_ls}).

ii) Next, suppose $w_1(\epsi_k)=\epsi_{m_1}>0$, $w_2(\epsi_k)=\epsi_{m_2}>0$. Note that $m_1>m_2\geq k$, because $w_1(\epsi_k)<w_2(\epsi_k)$ and $w_1(\epsi_i)=w_2(\epsi_i)$ for all $i$ from 1 to $k-1$. Here we put $k_2=m_1+1$, so $\wt w_1(k)\in C$ and $\wt w_2(k)\in B$. By Lemma~\ref{lemm:l_w_dash_A_n}~(i), $l''(w_1')\neq l''(w_2')$, where $w_i'=\wt w_is_{\eta_{k_1}-\eta_{k_2}}$, $i=1,2$. But $C_1'=C_2'$, a contradiction.

iii) Finally, suppose $w_1(\epsi_k)=-\epsi_{m_1}<0$, $w_2(\epsi_k)=-\epsi_{m_2}<0$. Note that $m_2>m_1>k$, because $w_1(\epsi_k)<w_2(\epsi_k)$ and $w_1(\epsi_i)=w_2(\epsi_i)$ for all $i$ from 1 to $k-1$. Here we put $k_2=m_2+1$, so $\wt w_1(k)\in B^-$ and $\wt w_2(k)\in C^-$. By Lemma~\ref{lemm:l_w_dash_A_n}~(ii), $l''(w_1')\neq l''(w_2')$, where $w_i'=\wt w_is_{\eta_{k_1}+\eta_{k_2}}$, $i=1,2$. On the other hand, $C_1'=C_2'$. This contradicts (\ref{formula:if_cones_then_ls}). The result follows.\hspace{\fill}$\square$\par

\nota{Actually, for $\Phi=B_n$, one can introduce the notion of basic involution literally as for $D_n$. It is easy to check that the previous proposition is true for basic involutions in type $B_n$.}

\end{document}